\documentclass[12pt, a4paper]{amsart}
\usepackage[margin=1in]{geometry}
\usepackage{mathtools}
\usepackage{amsthm}
\usepackage{amssymb}
\usepackage{amsrefs}
\usepackage{graphicx}
\usepackage{tikz}
	\usetikzlibrary{decorations.pathreplacing}
	\usetikzlibrary{patterns}
\usepackage{tikz-cd}
\usepackage{caption}
\usepackage{csquotes}
\usepackage[usestackEOL]{stackengine}
\usepackage{scalerel}
\usepackage{calc}
\usepackage{accents}

\usepackage{titlesec}
\titleformat*{\paragraph}{\scshape}
\titlespacing*{\paragraph}{0pt}{0pt}{0.5em}

\MakeOuterQuote{"}

\newcommand{\ra}{\rightarrow}
\newcommand{\pr}{\prime}
\newcommand{\de}{\partial}

\DeclareMathOperator{\ks}{ks}

\newtheoremstyle{thmstyle}		% name
    {\parskip}					% Space above
    {0pt}						% Space below
    {\slshape}					% Body font
    {}							% Indent amount
    {\scshape}					% Theorem head font
    {.}							% Punctuation after theorem head
    {.5em}						% Space after theorem head
    {}							% Theorem head spec (can be left empty, meaning ‘normal’)

\theoremstyle{thmstyle}

\newtheorem*{thm*}{Theorem}

\newtheorem*{defin*}{Definition}

\renewenvironment{proof}{{\noindent\slshape Proof.}}{\qed}

\makeatletter
\newsavebox\myboxA
\newsavebox\myboxB
\newlength\mylenA

\newcommand*\xoverline[2][0.75]{%
    \sbox{\myboxA}{$\m@th#2$}%
    \setbox\myboxB\null% Phantom box
    \ht\myboxB=\ht\myboxA%
    \dp\myboxB=\dp\myboxA%
    \wd\myboxB=#1\wd\myboxA% Scale phantom
    \sbox\myboxB{$\m@th\overline{\copy\myboxB}$}%  Overlined phantom
    \setlength\mylenA{\the\wd\myboxA}%   calc width diff
    \addtolength\mylenA{-\the\wd\myboxB}%
    \ifdim\wd\myboxB<\wd\myboxA%
       \rlap{\hskip 0.5\mylenA\usebox\myboxB}{\usebox\myboxA}%
    \else
        \hskip -0.5\mylenA\rlap{\usebox\myboxA}{\hskip 0.5\mylenA\usebox\myboxB}%
    \fi}
\makeatother

\newcommand{\CP}{{C\mkern-1.5muP}}
\newcommand{\CPbar}{\xoverline[0.75]{{C\mkern-1.5muP}}}
\newcommand{\cs}{\mathop{\#}}
\newcommand{\Z}{\mathbb Z}

\begin{document}

\title[Certifying a 4-manifold]{Certifying a Compact Topological 4-Manifold}

\author{Michael Freedman}
\address{\hskip-\parindent
	Michael Freedman\\
    Microsoft Research, Station Q, and Department of Mathematics\\
    University of California, Santa Barbara\\
    Santa Barbara, CA 93106\\}
\email{mfreedman@math.ucsb.edu}

\author{Daniele Zuddas}
\address{\hskip-\parindent
	Daniele Zuddas\\
	Lehrstuhl Mathematik VIII, Mathematisches Institut der Universit\"{a}t Bayreuth, NW II\\
	Universit\"{a}tsstr. 30, 95447 Bayreuth, Germany\\}
\email{zuddas@uni-bayreuth.de}

\begin{abstract}
We prove that compact topological 4-manifolds can be effectively presented by a finite amount of data.
\end{abstract}


\maketitle
The present paper is aimed at proving the following theorem.

\begin{thm*}
	Let $M$ be a compact topological 4-manifold. There is a finite data structure $C$ (a finite number of bits) which provably specifies $M$ up to homeomorphism. We call $C$ a "certificate" for $M$.
\end{thm*}

\paragraph{Background.}
Compact PL manifolds are manifestly specified by a finite amount of data. Up through dimension 3 the categories Top, PL, and Diff have the same objects so there is no issue of finding a certificate. For manifolds of dimension $\geq 5$ the surgery exact sequence reduces the specification of a structure $S(P)$ on a Poincar\'{e} space $P$ to homotopy theory (normal maps) and K-theoretic objects (the L-groups). A manifold structure (Top, PL, or Diff) can be certified in principle by giving a cell structure for the appropriate Poincar\'{e} space and an explicit solution for the homotopy and L-group problems. In dimension 4 PL and Diff have the same objects so the single open case is the question of a certificate for a compact topological 4-manifold.

\begin{proof}
	An early theorem of Lashof \cite{lashof} states that every noncompact connected 4-manifold admits a smooth structure. We will freely switch between smooth and PL language as we discuss smooth/PL manifolds of dimension 4 and 5. We use the language best suited to the construction at the moment; there are no essential differences in those dimensions.

	The broad outline of the proof is: (1) put a PL structure on $M^-$, $M$ minus a point, and (2) come up with enough smooth/PL data to prove the existence of a topologically flat 3-sphere cutting off the end of $M^-$. The certificate consists of the triangulation of a large, but finite, chunk of $M^-$, containing this $S^3$. Then, cutting along $S^3$ and gluing in $B^4$ recovers the homeomorphism type of $M$ uniquely.

	For the present, assume that $M$ is connected and orientable. To start we remove the Kirby-Siebenmann invariant $\ks$ which otherwise would cause us an inconvenience at the end of the proof. If $\ks(M) = 0$, let $M_-$ be the compact manifold obtained from $M$ by deleting the interior of a flat 4-ball $B^4 \subset$ int$(M)$. If $\ks(M) = 1$, let $M_-$ be the compact manifold obtained from $M$ by deleting the interior of $K$ from $M$, where $K$ is the compact contractible 4-manifold with boundary the Poincar\'{e} homology sphere. It is known that the double of $K$ is homeomorphic to $S^4$ so $K$ embeds in $B^4 \subset M$, thus we have such a copy of $K$ available to delete. In either case we have:
	\begin{equation}
		\ks(M_-) = 0
	\end{equation}

	Let $M^-_-$ be $M_- - \{\text{interior point}\}$. By Lashof's result \cite{lashof}, $M^-_-$ has a smooth structure extending the unique structure on its boundary; fix one. Let $D$ be its double along the entire boundary of $M^-_-$,
	\begin{equation}
		D = \operatorname{double}(M^-_-)
	\end{equation}

	Note that $D$ has two ends. We want to study (one of) these ends in the least distracting context so we do finitely many 1-surgeries to turn $D$ into a simply connected $D'$. Note that $D$ and $D'$ are smooth manifolds.

	To do this: first do 1-surgeries in both copies of $M^-_-$ to make them simply connected. If $M^-_-$ has $b$ boundary components $\pi_1(D)$ will now be a free group on $b-1$ generators. Kill this group by surgery on doubles of arcs joining each boundary component of $M^-_-$ to some reference boundary component. Now $D'$ is simply connected and may be assumed to still have the structure of a double. Thus, it admits an orientation-reversing involution, hence a $\Z_2$-action, which permutes the two halves. This extends over the end compactification $\xoverline D$ of $D'$ (which, by construction, is a topological manifold). In particular, $\xoverline D$ has zero signature and so its intersection form is indefinite or null.

	Therefore, the intersection form of $D'$ and of $\xoverline{D}$ is isomorphic to
	\begin{equation}
		\left (\bigoplus_k
		\begin{pmatrix}
			0 & 1 \\
			1 & 0 \\
		\end{pmatrix}\right ) \bigoplus \left (\bigoplus_j
		\begin{pmatrix}
			1 & 0 \\
			0 & -1 \\
		\end{pmatrix} \right).
	\end{equation}
	This follows by Serre's classification of integral unimodular indefinite symmetric bilinear forms \cite{serre}. Moreover, we can assume either $k=0$ or $j=0$.
	
	Since $\xoverline{D}$ is a double we have $\ks(\xoverline{D}) = 0$ but more importantly each half of $\xoverline{D}$, $\xoverline{D}/\Z_2$ also has vanishing $\ks$ invariant:
	\begin{equation}\label{kshalf/eqn}
		\ks(\xoverline{D}/\Z_2) = 0.
	\end{equation}

	By the classification of compact 1-connected topological 4-manifolds \cite{freedman} (and compactifying product ends and the inverse operation) we see that
	\begin{equation}\label{std/eqn}
		D' \underset{\text{top}}{\cong} \cs_k (S^2 \times S^2) \cs\, (\cs_j \CP^2) \cs\, (\cs_j \CPbar^2) - \{\text{2 pts}\}
	\end{equation}

	Denote the right hand side of (\ref{std/eqn}) by $E$; it is endowed with the canonical smooth structure. Observe that the smooth structure of $E$ clearly extends over the ends, while that of $D'$ may not be extendable.

	We may consider $D' \times [0,1]$ and ask if it has a smooth structure agreeing with $D' \times [0, \epsilon]$ at one end and $E \times (1-\epsilon, 1]$ at the other. The obstruction lies in
	\begin{equation}
		H^4 \left(D' \times [0,1], D' \times \{0, 1\}; \pi_3\!\left(\frac{\text{Top}}{\text{O}}\right)\right)\cong \pi_3\!\left(\frac{\text{Top}}{\text{O}}\right) \cong \Z_2
	\end{equation}

	The answer is: "Yes"\!. The obstruction precisely detects whether the corresponding ends of $D'$ and $E$ have the same Rochlin invariant \cite{freedman}. However all four Rochlin invariants vanish---this is why we took the $\ks$ invariant of $M$ into account when defining $M_-$ (compare with line (\ref{kshalf/eqn})). The fact we are using: if $X$ is a compact oriented topological spin 4-manifold of zero signature (possibly with boundary) and $X_S^-$ is a smoothing of $X$ minus an interior point, then:
	\begin{equation}\label{ks/eqn}
		\ks(X) = \operatorname{Rochlin}(\text{end } X_S^-)
	\end{equation}

	Indeed, both sides of (\ref{ks/eqn}) can be identified with the stable obstruction in $H^4(X, \partial X;\allowbreak \pi_3(\text{Top} / \text{PL})) \cong \Z_2$ to extend the smooth structure on $\partial X$ over $X$, see Kirby and Siebenmann \cite{KS77}. In our situation, it is enough to take $X \subset \xoverline{D}$ a small compact neighborhood of one end, which is bounded by a smooth connected 3-manifold in $D'$, to obtain $\text{ks}(X) = 0$.

	Let us denote such a smoothing of $D' \times [0,1]$ by $N$. Then, $N$ is a smooth proper $h$-cobordism from $D'$ to $E$. Fix a relative handle decomposition of $N$. As usual, we cancel all 0, 1, 4, and 5 handles of $N$. Everything about $N$ can be read off from the "middle level" $P$, a cross-section of $N$ above the 2-handles and below the 3-handles.

	Inside $M$ we see two collections of disjoint, locally finite, smooth spheres $A$ and $B$. $A$ consists of the descending (attaching) 2-spheres of the 3-handles and $B$ the ascending (belt) 2-spheres of the 2-handles. Since $N$ is a proper $h$-cobordism, algebraically we may assume that $A$ and $B$ intersect as follows:
	\begin{equation}
		A_i \cdot B_j = \delta_{ij}
	\end{equation}

	and also by a lemma of Casson \cite{casson}, up to finger moves we may assume that
	\begin{equation}
		\pi_1(P -(A \cup B)) \cong 0
	\end{equation}

	But geometrically there are, in addition to the desired intersection points, additional points which may be arranged in cancelling pairs, each pair on a Whitney circle, and all Whitney circles bounding locally flat 2-disks contained in disjoint smooth 6-stage Casson towers in $P$ which meet $A \cup B$ only along the Whitney circle.

	\begin{defin*}
		Fixing a gradient-like flow associated with a handle decomposition in any category (Top, PL, or Diff), a \textsl{generalized flow line} is a minimal closed invariant subset of the bidirectional flow (see Figure 1).
	\end{defin*}

	Now consider a smooth 3-sphere $S^3$ cutting off one of the ends of $E$ at the top of the $h$-cobordism. $S^3$ should be indeed very near an end, in that we do not want any of the "generalized flow" lines leaving $S^3$ to arrive in a portion of $D'$ modified by the above finitely many 1-surgeries; all such generalized flow lines should end in a part of $D'$ identified with a portion of (one copy of) the original manifold $M$. Similarly, all generalized flow lines through $S^3$ should be disjoint from some fixed compact 2-complex $T$ in $D'$ carrying $H_2(D'; \Z)$. Notice that such a 3-sphere exists topologically, but not necessarily smoothly, also in $D'$. Our goal is to represent such topologically locally flat 3-sphere in $D'$ starting from the smooth sphere $S^3 \subset E$.

	\begin{figure}[ht]
		\centering
		\begin{tikzpicture}[scale = 1.2]
			\draw  (-0.1,0) rectangle (2.4,0.2);
			\draw (2.4,0.2) -- (3.7,1.6);
			\draw (2.4,0) -- (3.7,1.4);
			\draw (3.7,1.4) -- (3.7,1.6);

			\draw (5.9, 0.5) to [out = 85, in = 225] (6.3, 2.1) to [out = 35, in = 180] (6.6,2.2) to [out = 0, in = 145] (6.9, 2.1) to [out = -45, in = 95] (7.3, 0.5);
			\draw (7.3, 0.5) to [out = 25, in = -60] (7.8, 1.4);
			\draw [dashed] (7.8, 1.4) to [out = 120, in = 70] (5.7, 0.7);
			\draw (7.85, 1.3) to [out = 100, in = 0] (6.8,2.7) to [out = 180, in = 105] (5.7, 0.7) to [out = 295, in = 160] (5.9, 0.5);

			\fill (6.6, 2.2) circle[radius = 1.5pt];
			\fill (7,2.7) circle[radius = 1.5pt];
			\draw (6.6, 2.2) -- (7,2.7);
			\draw [pattern = vertical lines] (7,3) to [out = -135, in = 90] (6.4, 2) to [out = -90, in = 225] (6.8,2.1) to [out = 45, in = -90] (7,2.4) -- (7,3);

			\draw [gray, thin] (7.5, 0.62) to [out = 90, in = -30] (6.95, 2.25);
			\draw [gray, thin] (7.75, 0.86) to [out = 90, in = -30] (7, 2.45);
			\draw [gray, thin] (6.45, 2.3) to [out = 190, in = 75] (5.65, 1.3);
			\draw [gray, thin] (6.55, 2.5) to [out = 185, in = 55] (5.87, 2.1);
			\draw [gray, thin] (6.45, 1.95) to [out = 230, in = 80] (6.25, 1.4);
			\draw [gray, thin] (6.8, 2.1) to [out = 230, in = 80] (6.55, 1.56);
			\node at (4.6,1.7) {\resizebox{1.6em}{!}{$\hookleftarrow$}};

			\draw (0.9, 0.5) to [out = 85, in = 225] (1.3, 2.1) to [out = 35, in = 180] (1.6,2.2) to [out = 0, in = 145] (1.9, 2.1) to [out = -45, in = 95] (2.3, 0.5);
			\draw (2.3, 0.5) to [out = 25, in = -60] (2.8, 1.4);
			\draw [dashed] (2.8, 1.4) to [out = 120, in = 70] (0.7, 0.7);
			\draw (2.85, 1.3) to [out = 100, in = 0] (1.8,2.7) to [out = 180, in = 105] (0.7, 0.7) to [out = 295, in = 160] (0.9, 0.5);
			\draw  (0.9, 0.5) ellipse (0.15 and 0.1);
			\draw  (2.3, 0.5) ellipse (0.15 and 0.1);
			\draw (0.75, 0.5) to [out = 90, in = 180] (1.6, 2.35) to [out = 0, in = 90] (2.45, 0.5);
			\draw (1.05, 0.5) to [out = 90, in = 265] (1.1,1.2) to [out = 85, in = 180] (1.6, 2.05) to [out = 0, in = 95] (2.1, 1.2) to [out = 275, in = 90] (2.15, 0.5);

			\draw (3.7, 1.6) -- (2.8, 1.6);
			\draw (-0.1,0.2) -- (0.66,0.9);
			\node at (9.4,1.8) {generalized};
			\node at (9.4,1.4) {flow line};
			\draw [->] (0,1.4) -- (0.6,1.7);
			\draw [->] (0,0.8) -- (0.9,1.3);
			\node at (-0.8,1.3) {2-handle};
			\node at (-0.8,0.7) {1-handle};
		\end{tikzpicture}
		\caption{}
	\end{figure}

	We now identify a crucial subset $Y \subset P$ of the middle level. Let $A'$ and $B'$ be those spheres within $A$ and $B$ respectively whose generalized flow lines meet $S^3 \subset E$. If $A_i \in A'$ it is easy to see that $B_i \in B'$, however $B_i \in B'$ does not imply $A_i \in A'$ (the asymmetry comes from the fact $S^3$ is at the top of the $h$-cobordism). To correct this imbalance let $A^{\pr\pr}$ be $A'$ union all spheres in $A$ which are duals of spheres in $B'$, i.e. $B_j \in B' \Rightarrow A_j \in A^{\pr\pr}$. Let $C$ be the union of all 6-stage Casson towers on all Whitney circles meeting $A^{\pr\pr} \cup B'$. Note the Whitney circle might pair double points between $A^{\pr\pr} \cup B'$ and $A \cup B - (A^{\pr\pr} \cup B')$. Now $Y$ is defined as a closed regular neighborhood

	\[
	Y = \mathcal{N}(A^{\pr\pr}\cup B' \cup C) \subset P
	\]

	It is a key result \cite{freedman} that every 6-stage Casson tower, a combinatorially explicit object, contains a topological 2-handle $H^2_{\text{top}}$, which is not an explicit finitistic object. $Y$ can be described with finite data but we are about to use the topological Whitney disk $W \subset H^2_{\text{top}}$, which is the core of $H^2_{\text{top}}$, to construct a flat 3-sphere cutting off an end of $D'$ and similarly cutting off the end of $M^-_-$.

	Up to this point all flow lines are smooth or PL. Now perform an ambient topological isotopy of $B'$ within $Y$ using the Whitney disks in the topological 2-handles within $C$ to make $A^{\pr\pr}$ and $B'_{\text{top}}$, the isotoped $B'$, have geometric $\delta_{ij}$ intersection:
	\begin{equation}
		|A_i \cap (B_{\text{top}})_j| = \delta_{ij}
	\end{equation}
	for all $A_i \in A^{\pr\pr}$ and $(B_{\text{top}})_j \in B'_\text{top}$.

	After this isotopy the generalized flow lines in $N$ are now topological, not smooth, but now all generalized flow lines starting in $S^3 \subset E$ are either intervals running from $E$ to $D'$ or else meet $E$ in topologically flat 2-disks. This statement is the essence of the famous Morse cancellation lemma, and sketched in two lower dimension, in Figure 1. Said another way, an ascending sphere is punctured by its encounter with its dual (descending) sphere and arrives at $E$, the top of $N$, as a flat 2-disk.

	Our goal is that $S^3$ should meet only ordinary (interval) flow lines. But since the generalized flow lines it meets intersect $E$ in finitely many disjoint flat 2-disks it is easy to find an ambient topological isotopy of $S^3$ to $S'$ where $S'$ now only meets ordinary interval flow lines. The union of these ordinary flow lines is a topologically flat $(S^3 \times [0,1]; S', S'') \subset (N; E, D')$. This constructs a topologically flat 3-sphere $S''$ far out an end of $D'$ and cutting off that end.

	By placing $S^3$ sufficiently near an end of $E$ we have ensured that the region in $D'$ related to it through generalized flow lines is nearer the end than both the finitely many $S^2 \times D^2$'s added during surgery $D \ra D'$, and the 2-complex $T$ carrying $H_2(D'; \Z)$. This ensures that any topologically flat 3-sphere in $D'$ lying within the generalized flow lines through $Y$ and which cuts off an end, cuts off a punctured contractible manifold from $D'$. Homeomorphically this contractible manifold, since its boundary is a sphere, must be a 4-ball $B^4$ by the 4-dimensional topological Poincar\'{e} conjecture. Moreover, $S''$ can be identified with a locally flat sphere in $M^-_-$, still denoted by $S''$, that cuts off the end.

	\paragraph{Claim.} The manifold $M$ can be uniquely recovered up to homeomorphism as:
	\begin{equation}
		M \underset{\text{top}}{\cong} M^\ast \cup (B^4 \text{ or } K) \cup B^4
	\end{equation}
	where $M^\ast$ is the compact component of $M^-_-$ cut open along $S''$, the first "$B^4$ or $K$" restores the bit removed at the beginning to ensure $\ks(M_-) = 0$ and the final $B^4$ caps off $S'$.

	The claim follows from two applications of the fact \cite{freedman} that if $J$ is a compact contractible manifold then every homeomorphism of $\de J$ extends to a homeomorphism of $J$, and moreover if $J'$ is another contractable manifold such that $\de J \cong \de J'$ then any homeomorphism $\partial J \to \partial J'$ extends to a homeomorphism $J \ra J'$. It is applied once to regluing $B^4$ or $K$ and once to capping off the severed end of $M^-_-$ with $B^4$, ensuring that the result is independent of the gluing homeomorphism.

	The finite "certificate" for the homeomorphism type of $M$ consists of the following pieces of information:
	\begin{enumerate}
		\item The initial binary choice: removal of $B^4$ versus $K$ from a ball in int($M$) to obtain $M_-$
		\item A triangulation of the complement of a small neighborhood of infinity of the punctured $M^-_-$. The neighborhood should be sufficiently small that a collection of 1 and 2-cells presenting $\pi_1(M)$ and generating a basis for the second homology $H_2(M; \Z)$ should exist outside the neighborhood of infinity.
		\item The portion of the smooth proper $h$-cobordism $N$ consisting of the generalized flow lines through a PL neighborhood of $Y$ and of $S^3 \subset E$. This is a large enough piece of $N$ to prove that the sphere $S^3$ in $E$ can be isotoped to $S'$ and then pushed down to $M^-_-$ inside $D'$.
	\end{enumerate}
	This information suffices to reconstruct $M$ up to homeomorphisms.

	This completes the proof except for the nonorientable case. Assume $M$ is nonorientable and that the first Stiefel-Whitney class is classified by a map $f$:
$$\begin{tikzcd}
	M  \ar{r}{f} & RP^\infty \\
	L \ar[r]  \ar[u, hook] & RP^{\infty-1}. \ar[u, hook]
\end{tikzcd}$$

By topological transversality \cite{fq} we may assume $f$ is transverse to $RP^{\infty - 1}$ with preimage a submanifold $L$. Let $Q$ be $M$ cut open along $L$. Then $Q$ is orientable and by the previous argument admits a certificate {(1), (2), and (3)} as above. Now to this add a fourth term: 4) a description of the fixed point free involution $\iota$ of $\widetilde{L}$ so that $\widetilde{L}/\iota = L$, where $\widetilde{L}$ is the boundary of a tubular neighborhood of $L$ in $M$. This is possible since $L$ is a 3-manifold and hence can be explicitly triangulated. Then $M$ can be reconstructed as
	$$M \underset{\text{top}}{\cong} Q/\iota,$$ where the identifications by $\iota$ are within $\widetilde L$.
\end{proof}

\newpage

\section*{Acknowledgments}
The second author acknowledges support of the 2013 ERC Advanced Research Grant 340258 TADMICAMT. He is also a member of GNSAGA -- Istituto Nazionale di Alta Matematica ``Francesco Severi'', Italy.

%\newpage

\end{document}